\documentclass[letterpaper,10pt,conference]{IEEEtran}
\IEEEoverridecommandlockouts
\usepackage{cite}
\usepackage{amsmath,amssymb,amsfonts}
\usepackage{algorithm}
\usepackage{algorithmic}
\usepackage{graphicx}
\usepackage{textcomp}
\usepackage{xcolor}
\usepackage{epstopdf}

\usepackage{atbegshi}
\def\BibTeX{{\rm B\kern-.05em{\sc i\kern-.025em b}\kern-.08em
    T\kern-.1667em\lower.7ex\hbox{E}\kern-.125emX}}
\usepackage[letterpaper, left=19.2mm, right=19.1mm, top=20.2mm, bottom=20.5mm]{geometry}
\begin{document}

\allowdisplaybreaks[4]

\AtBeginShipout{
    \ifnum\value{page}=1
        \newgeometry{left=19.2mm, right=19.1mm, top=25.2mm, bottom=20.5mm}
    \fi 
}
\AtBeginShipout{
    \ifnum\value{page}=2
        \newgeometry{left=19.2mm, right=19.1mm, top=20.2mm, bottom=20.5mm}
    \fi 
}

\title{\LARGE \bf Deep Learning Model Predictive Control for Deep Brain Stimulation in Parkinson's Disease\\

\author{Sebastian Steffen$^\ast$ \thanks{$^\ast$Department of Engineering Science, University of Oxford, OX1 3PJ, UK. {\tt\footnotesize \{sebastian.steffen,mark.cannon\}@eng.ox.ac.uk}}
\and Mark Cannon$^\ast$
}

\thanks{This work was supported by the UKRI Engineering and Physical Sciences Research Council}
}
\maketitle

\begin{abstract}
We present a nonlinear data-driven Model Predictive Control (MPC) algorithm for deep brain stimulation (DBS) for the treatment of Parkinson's disease (PD). Although DBS is typically implemented in open-loop, closed-loop DBS (CLDBS) uses the amplitude of neural oscillations in specific frequency bands (e.g.~beta 13-30\,Hz) as a feedback signal, resulting in improved treatment outcomes with reduced side effects and slower rates of patient habituation to stimulation.
To date, CLDBS has only been implemented \textit{in vivo} with simple algorithms such as proportional, proportional-integral, and thresholded switching control. 
Our approach employs a multi-step predictor based on differences of input-convex neural networks to model the future evolution of beta oscillations. The use of a multi-step predictor enhances prediction accuracy over the optimization horizon and simplifies online computation. In tests using a simulated model of beta-band activity response and data from PD patients, we achieve reductions of more than 20\% in both tracking error and control activity in comparison with existing CLDBS algorithms. The proposed control strategy provides a generalizable data-driven technique that can be applied to the treatment of PD and other diseases targeted by CLDBS, as well as to other neuromodulation techniques.
\end{abstract}

\begin{IEEEkeywords}
neuromodulation, data-driven control, model predictive control, nonlinear systems
\end{IEEEkeywords}

\section{Introduction}
Deep Brain Stimulation (DBS) is a treatment for various neurological and psychiatric disorders that involves the surgical implantation of electrodes into specific structures deep within the brain.
DBS devices deliver pulses of electric current to disrupt pathological activity in the central nervous system. The technique is currently used to treat essential tremor, Parkinson's disease (PD), and epilepsy~\cite{sarica:2023}. It is also being trialed for treatment-resistant depression and obsessive-compulsive disorder, among other conditions~\cite{groppa:2024,Sellers:2024}. Typically, DBS operates in an open-loop mode with a fixed stimulation pattern of constant amplitude, frequency and pulse-width, which results in high stimulation levels compared to targeted stimulation algorithms~\cite{schiff:2010}. Excessive stimulation can increase side-effects and lead to a more rapid decline in the efficacy of treatment due to habituation~\cite{fleming:2020}. This has motivated extensive research into closed-loop DBS (CLDBS), leading to recent approval of the technique for PD patients\cite{medtronic:2025,FDA:2025}. 

Typically, CLDBS algorithms typically modulate stimulation amplitude as a function of disease biomarkers measured by the implant\cite{Acharya:2025}. Symptoms of PD are associated with bursts in amplitude of so-called beta-band oscillations (13-30Hz) in population-level neural activity known as local field potentials (LFP)~\cite{silberstein:2005,tinkerhauser:2020}. To date, \textit{in vivo} CLDBS has been limited to simple algorithms such as on/off switching control \cite{he:2023,arlotti:2018}, proportional or proportional-integral (PI) control \cite{schmidt:2023}, or dual-threshold control \cite{velisar:2019,Stanslaski:2024}.

Simulation-based studies allow a wider range of control strategies to be investigated. Some approaches incorporate additional information, such as muscle activity, to select between PI controllers tuned for different operating points~\cite{fleming:2023}. Others employ different control objectives, e.g., minimizing the duration of periods of high pathological activity~\cite{petrucci:2020}. Another research direction is to use model-based feedback control laws. In~\cite{yang2018control} an LQR controller is designed to control the frequency of DBS pulses based on a model derived from input-output data, and in~\cite{Fang:2024} this is augmented with nonlinear terms to improve model accuracy, and the controller contains additional terms to cancel the nonlinearity. 

Several studies have proposed optimal predictive control schemes, although no consensus has emerged on the type of model to use. The majority of schemes fit a linear model to the dynamics of the relevant biomarker, for example, ARX models~\cite{haddock:2017} or state-space models fit using subspace methods~\cite{ahmadipour:2021}. In our previous work~\cite{steffen:2024}, we proposed an augmented model comprising an online-identified linear model of biomarker activity and a second-order model of the response to stimulation based on averaged patient data. While linear models can provide a simpler formulation of the optimal control problem, they do not accurately capture the inherently nonlinear dynamics of neural systems \cite{Rubin:2004}. Very few studies have investigated nonlinear predictive control. In \cite{su:2018}, a Volterra series is used to represent the nonlinear dynamics of the patient response to DBS. To the best of our knowledge, no studies have considered the use of more expressive nonlinear models such as 
feedforward or recurrent 
neural networks in the context of CLDBS. 

If the nonlinear model can be represented as a difference of convex (DC) functions, the online MPC optimization problem can be solved efficiently as a sequence of convex problems~\cite{DoffSotta:2022,Lishkova:2025}.
These are derived via partial linearization of the system model around nominal predicted trajectories.  
The DC model representation provides tight bounds on linearization errors and the nominal trajectories are successively updated using the most recent solution estimate. For example, \cite{krausch:2025deeplearningadaptivemodel} uses input-convex neural networks (ICNN~\cite{amos:2017}) to apply robust MPC to a batch bioreactor. The resulting DC-MPC controller, like other robust nonlinear MPC approaches, requires the propagation of linearization errors over the prediction horizon to ensure satisfaction of constraints. As discussed in~\cite{koehler:2022}, a multi-step predictor model simplifies the construction of robust tubes bounding the future model state.

\subsubsection*{Contribution}
This paper provides the following contributions to the state of the art in CLDBS. 
\begin{itemize}
\item We formulate a nonlinear optimal control problem for CLDBS subject to input and output constraints, using a multi-step predictor model defined for each prediction step by the difference of a pair of ICNNs.
\item We solve this problem using sequential convex programming with tight bounds on approximation errors.
\item We construct the control law as a robust tube MPC strategy that explicitly accounts for linearization errors and uncertainty in the predictor model.
\item The use of a multi-step predictor simplifies online computation by avoiding the need for recursive propagatation of linearization errors over the prediction horizon. 
\end{itemize}

The remainder of this paper is structured as follows. 
Section~\ref{sec:problem} outlines the problem formulation for CLDBS. Section~\ref{sec:controller} describes the proposed control law. Section~\ref{sec:simulations} describes numerical simulations and provides a discussion of their results. Finally, Section~\ref{sec:conclusions} provides concluding remarks and some perspectives on future work.

\section{Problem Formulation} \label{sec:problem}
The envelope of the beta-band oscillations measured at the implant site in the subthalamic nucleus at time $t$ is denoted $y(t) \in \mathbb{R}^{+}$. We assume that the dynamics of the biomarker $y(t)$ are given by an unknown, nonlinear (and possibly time-varying) dynamical system of the form
\begin{equation}\label{eq:cont_time_system}
    \dot{y}(t) = f(y(t), u(t), t)
\end{equation}
where the control input $u(t)$ is the product of the applied stimulation amplitude in volts, pulse width in seconds and frequency in Hertz, and is constrained to lie in the range $u \in[0, u_{\mathrm{max}}]$ which is determined for each patient by a clinician. Furthermore, the rate of change of input is also constrained, lying in the range $\dot{u}(t) \in [-\dot{u}_{\mathrm{max}},  \dot{u}_{\mathrm{max}}]$. We assume that $f$ is a continuous, twice-differentiable function, and can therefore be represented with arbitrary accuracy by a difference of convex functions. 

The goal of CLDBS is to suppress beta-band activity exceeding a pathological level, which we denote $y_{0}$, while minimizing the stimulation energy. This suggests a cost index for the optimal control law defined over a horizon of length $T$ of the form:
\begin{gather}
\label{eqn:cost_fn}
\int_0^T \Bigl( \phi([y(t)-y_0]_{\geq 0}) + R u(t)^{2} \Bigr) \, dt
\\
\label{eqn:positive_cost}
    [y(t) - y_0]_{\geq 0} = 
    \begin{cases}
        y(t) - y_{0}, & \text{if}\ y(t) \geq y_{0} \\
        0 & \text{otherwise}
    \end{cases}
\end{gather}
where $\phi:\mathbb{R}_{\geq 0} \to \mathbb{R}_{\geq 0}$ is a monotonically non-decreasing function, and $R$ is a  positive control weighting.

\section{Control Law} \label{sec:controller}
\subsection{Difference of Convex Functions Neural Network Model}
To derive an MPC strategy we construct a discrete time model of the system~\eqref{eq:cont_time_system}. 
For the predictor model in discrete time, we define a separate neural network for each of the $N$ steps of the prediction horizon:
\begin{equation} \label{eq:nn_dynamics}
    \begin{bmatrix}
        y_{k+1} \\ 
        y_{k+2} \\ 
        \vdots \\ 
        y_{k+N}
    \end{bmatrix} = 
    \begin{bmatrix}
        f_{1}(z_{k}, u_{k}) \\
        f_{2}(z_{k}, u_{k:k+1}) \\ 
        \vdots \\
        f_{N}(z_{k}, u_{k: k+N-1})
    \end{bmatrix} +
    \begin{bmatrix}
        w_{1} \\
        w_{2} \\
        \vdots \\
        w_{N}
    \end{bmatrix}
\end{equation} 
where $k$ is the discrete time index for a given sampling interval and each function $f_i$ is a difference of convex functions 
\begin{equation}
f_{i}(z_{k}, u_{k:k+i-1}) = f_{i,1}(z_{k}, u_{k:k+i-1}) - f_{i,2}(z_{k}, u_{k:k+i-1}) .
\end{equation}
The arguments of $f_i,f_{i,1},f_{i,2}$ are $z_{k} = [y_{\mathrm{past}, k} \  u_{\mathrm{past},k}]^{T}$ which contains $n_{y}$ past observations of the system output and $n_{u}$ control inputs, i.e.~$y_{\mathrm{past},k} = [
    y_{k} \ y_{k-1} \ \cdots \ y_{k-n_{y}+1} ]$ 
and $u_{\mathrm{past},k} = [u_{k-1} \ u_{k-2} \ \cdots \ u_{k-n_{u}}]$, and the sequence of control inputs from timestep $k$ until timestep $k+i-1$, i.e.~$u_{k:k+i-1} = [u_{k} \ u_{k+1} \ \cdots \ u_{k+i-1}]^{T}$, and the disturbance $w_{i} \in \mathbb{W}_{i} = [w_{i, \mathrm{min}},\ w_{i,\mathrm{max}}]$ accounts for modeling error.

The functions $f_{i,1},f_{i,2}$ are each specified by ICNNs with the structure proposed in \cite{amos:2017}, consisting of a series of fully connected hidden layers, with additional skip connections from the input to the hidden layers. The use of rectified linear units as activation functions, together with non-negative constraints on the dense weights between hidden layers ensures convexity with respect to the network inputs. 

\subsection{DCNN Tube MPC}
We propose a robust MPC law obtained by minimizing a discrete-time, convex approximations of the cost~\eqref{eqn:cost_fn}. Each convex approximation is obtained by linearizing the concave parts of the model~\eqref{eq:nn_dynamics} around a nominal trajectory $(\vec{y}^{\,0}_{k}, \vec{u}^{\,0}_{k})$, where $\vec{y}^{\,0}_{k} = [
    y^0_{k+1} \ \cdots \ y^0_{k+N} ]$, $\vec{u}^{\,0}_{k} = [
    u^0_{k} \ \cdots \ u^0_{k+N-1} ]$.
We define perturbations, $s_{k} = y_{k} - y^0_{k}$, and $v_{k} = u_{k} - u^0_{k}$, and sets $\mathbb{S}_{k+i} = [s_{k+i,\mathrm{min}}, s_{k+i,\mathrm{max}}]$ bounding the predicted deviation of $y_{k+i}$ from $y^0_{k+i}$ at time $k+i$. The DC property of $f_i$ allows us to find tight bounds on the perturbations satisfying: 
\begin{align} 
s_{k+i, \mathrm{max}} &\geq f_{i,1}(z_{k}, u_{k:k + i-1}) - f_{i,2}(z_{k},u^0_{k:k+i-1}) \nonumber\\
&\quad - f_{i,2}^\prime (z_k, u^0_{k:k+i-1}) v_{k:k+i-1} + \underset{w\in \mathbb{W}_{i}}{\mathrm{max}}w - y^0_{k+i}
\label{eq:upper_bound}
\\
s_{k+i, \mathrm{min}} &\leq {-f_{i,2}}(z_{k},u_{k:k+i-1}) + f_{i,1}(z_{k}, u^0_{k:k + i-1})  \nonumber\\
&\quad + f_{i,1}^\prime(z_k, u^0_{k:k+i-1}) v_{k:k+i-1} + \underset{w\in \mathbb{W}_{i}}{\mathrm{min}}w - y^0_{k+i} 
\label{eq:lower_bound}
\end{align}
where $f_{i,1}^\prime(z_k,u^0_{k:k+i-1})$ and $f_{i,2}^\prime(z_k,u^0_{k:k+i-1})$ denote the Jacobian matrices $\partial f_{i,1}/\partial u$ and $\partial f_{i,2}/\partial u$ evaluated along $\vec{u}^{\,0}_k$.

\subsection{Optimal MPC Problem}
We define the worst-case MPC cost
\begin{equation} \label{eqn:dcnn_cost}
    J(z_{k}, \vec{u}_{k}, \vec{\mathbb{S}}_k) = \sum_{i=0}^{N} 
    \!\max_{s_{k+i}\in \mathbb{S}_{k+i}} \! 
    \Bigl( Q[y^0_{k+i} + s_{k+i}- \beta_{0}]^{2}_{\geq0} + Ru_{k+i}^{2} \Bigr)
\end{equation}
where $Q$ and $R$ are positive, scalar weights, and 
$[\cdot]_{\geq 0}$ denotes projection onto the positive orthant so that the term $Q[y^0_{k+i} + s_{k+i}- \beta_{0}]^{2}_{\geq0}$ only penalizes the predicted values of $y$ exceeding the threshold $\beta_{0}$.
We denote $\vec{u}^{\,*}_{k}$ as the solution of the convex program:
\begin{equation} \label{eq:optimal_u}
    \underset{\vec{u}_k, \vec{\mathbb{S}}_k}{\mathrm{minimize}}\hspace{0.2cm} J(z_{k},\vec{u}_{k},\vec{\mathbb{S}}_k)
\end{equation} subject to \eqref{eq:upper_bound}, \eqref{eq:lower_bound} and
\begin{align*} 
    y^0_{k+i} + \mathbb{S}_{k+i} &\subseteq \mathbb{Y} = [y_{\mathrm{min}}, y_{\mathrm{max}}]
\\
    u_{k+i} &\in \mathbb{U} = [u_{\mathrm{min}}, u_{\mathrm{max}}]
\\
    \Delta u_{k+i} &\in \Delta \mathbb{U} = [-\Delta u_{\mathrm{max}}, \Delta u_{\mathrm{max}}]
\end{align*} 
for all $i \in \{0, \dots, N-1\}$. Here $\mathbb{Y}$, $\mathbb{U}$ and $\Delta\mathbb{U}$ are the output, control input, and input rate constraint sets, with $\Delta u_{k} = u_{k} - u_{k-1}$.
The solution of this optimal control problem is used to update $\vec{u}^{\, 0}_{k}$ and hence compute the nominal predicted trajectory $\vec{y}^{\, 0}_{k}$ using \eqref{eq:nn_dynamics} with $w_i=0$ for all $i$. We iteratively update this solution until convergence (indicated by the change in the optimal cost falling below a threshold $\Delta J_{min}$), or until the maximum number of iterations ($\mathit{maxiters}$) is reached, as outlined in Algorithm \ref{alg:dcnn_tmpc}. 

\begin{algorithm}
\caption{Multi-step DCNN TMPC}
\label{alg:dcnn_tmpc}
\begin{algorithmic}[1]  
  \REQUIRE $z_{k}$, feasible $\vec{y}^{\, 0}_{k}$ and $\vec{u}^{\, 0}_{k}$ 
  \STATE {$j \gets 1$, $\Delta J \gets 10^6$, $J_0 \gets 10^6$ }
  \WHILE{$j \leq \mathit{maxiters}$ \textbf{and} $\Delta J > \Delta J_{min}$}
    \STATE Evaluate the Jacobian matrices $f_{i,1}^\prime$ and $f_{i,2}^\prime$ using $\vec{u}^{\,0}_k$
    \STATE Solve problem~\eqref{eq:optimal_u} for $\vec{u}^{\, *}_{k, j}$, given $\vec{y}_k^{\, 0}$ and $\vec{u}_k^{\, 0}$
    \STATE $J_{j} \gets J(z_{k}, \vec{u}^{\, *}_{k, j})$, $\Delta J \gets J_{j} - J_{j-1}$\\
    $\vec{u}^{\, 0}_{k}\gets \vec{u}^{\, *}_{k, j}$, $y^{0}_{k+i} \gets f_{i}(z_{k}, u^{0}_{k:k+i-1})$, $i=1,\ldots,N$ \\
    \STATE $j \gets j+1$
  \ENDWHILE  
  \RETURN $\vec{u}^{\, *}_{k}$, $\vec{y}^{\, 0}_{k}$, $\vec{u}^{\,0}_{k}$
\end{algorithmic}
\end{algorithm}

\section{Numerical Simulations} \label{sec:simulations}
We evaluate the performance of the proposed control scheme using LFP data from four individual Parkinsonian patients who underwent DBS surgery at either King's College Hospital or St George's Hospital in London. The LFP data was gathered while the patients were in a resting state, with stimulation switched off, for periods ranging from $15$ to $30$ minutes. The noisy LFP data was band-pass filtered between 18 and 24 Hz using a 6th order causal Butterworth filter, and the envelope was extracted via continuous wavelet transform. As in~\cite{steffen:2024}, we use here a synthetic model of the patients' beta-band activity in response to applied stimulation. The envelope of the beta-band oscillations is related to the nominal activity $y_{\beta}$ (i.e. the brain's activity in absence of any stimulation) and the DBS attenuation effect $\eta(t)$ as follows,
\begin{equation} \label{eq:beta_response}
    y(t) = y_{\beta(t)} \cdot e^{-\eta(u(t))}
\end{equation}
The stimulation response $\eta(t)$ is represented by a second-order continuous-time system,
\begin{equation}\label{eqn:stim_dynamics}
\begin{aligned}
    \dot{x}_{c}(t) &= 
    \begin{bmatrix} -1/\tau_{1}(t) & 0 \\
                    g(t)/\tau_{2}(t) & -1/\tau_{2}(t)
    \end{bmatrix}
    x_{c}(t) + \begin{bmatrix} g(t)/\tau_{1}(t) \\ 0\end{bmatrix} u(t) \\
    \eta(t) &= \begin{bmatrix}
        0 & 1
    \end{bmatrix} x_{c}(t)
\end{aligned}
\end{equation}
where the parameters $\tau_1(t)$, $\tau_2(t)$ and $g(t)$ differ across patients, and also vary over time due to short-term changes in the patient's activity and long-term changes due to disease progression. The average values for these parameters are
\[
\bar{g} = 62.11, \quad 
\bar{\tau}_{1} = 0.05, \quad
\bar{\tau}_{2} = 0.25.
\]

\subsection{Stimulation response and model training} \label{ssec:stim_train}
Simulations use a discrete-time representation of the synthetic beta response model (\ref{eqn:stim_dynamics}) employing a zero-order hold with sample rate $50$\,Hz. The parameters $\tau_{1}(t)$, $\tau_{2}(t)$ and $g(t)$ of the  model~(\ref{eqn:stim_dynamics}) vary in a random walk, with steps at each sampling instant of no more than 2.5\% of their nominal value ($\bar{\tau}_1$, $\bar{\tau}_2$ or $\bar{g}$) and with the constraint that the total variation does not exceed 40\% of the nominal value. 

The ICNNs were implemented using Keras \cite{chollet:2015} and trained on synthetically modulated LFP data sampled at 50 Hz. Training inputs were constructed using a psuedo-random binary sequence alternating between $-\Delta u_{max}$ and $\Delta u_{max}$ to define the incremental signal $\Delta u_{k} = u_{k} - u_{k-1}$ subject to $u_k\in [0, u_{max}]$. Networks were initially trained on $10^5$ samples of synthetically modulated trajectories of LFP data taken from three patients (to create the `pre-trained' model), then further trained on $3.2\times 10^4$ samples from the fourth patient (to create the `refined' model), with $10^4$ samples kept as the test set and with an offset of $3000$ samples between test and training sets. In all experiments, we compare the performance of the pre-trained model and the refined model to test the ability of the model to generalize to beta activity of an unseen patient. 100 epochs of training 5 models on an Apple M1 Pro CPU took 9 minutes, 20 seconds. 

\subsection{Multi-step Prediction Accuracy}
\begin{figure}
    \centering
    \includegraphics[width=0.85\linewidth,trim={0 3.5mm 0 3mm},clip]{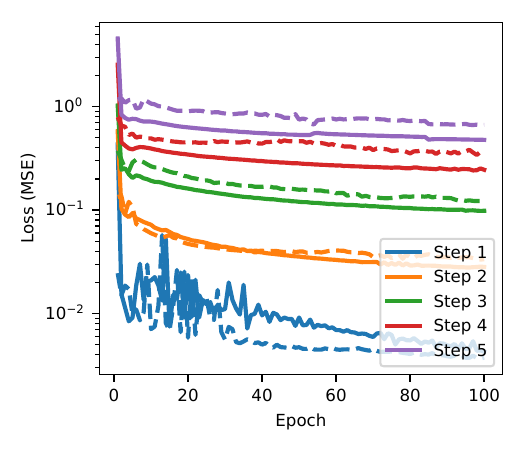}
    \caption{Training (solid) and validation (dashed) loss against training epoch for the multi-step predictor model over 1-5 steps}
    \label{fig:dcnn_training}
\end{figure}

\begin{figure}
    \centering
    \includegraphics[width=0.85\linewidth,trim={0 3.5mm 0 3mm},clip]{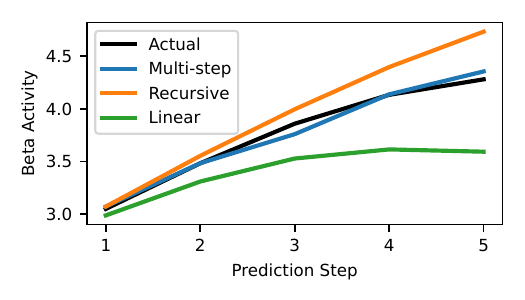}
    \caption{Example of predicted trajectories over a 5-step horizon for the multi-step DCNN, and recursive predictions with the single-step DCNN and linear ARI model, compared to the true beta activity}
    \label{fig:predicted_beta_trajectories}
\end{figure}

\begin{figure}
    \centering
    \includegraphics[width=0.85\linewidth,trim={0 3.5mm 0 3mm},clip]{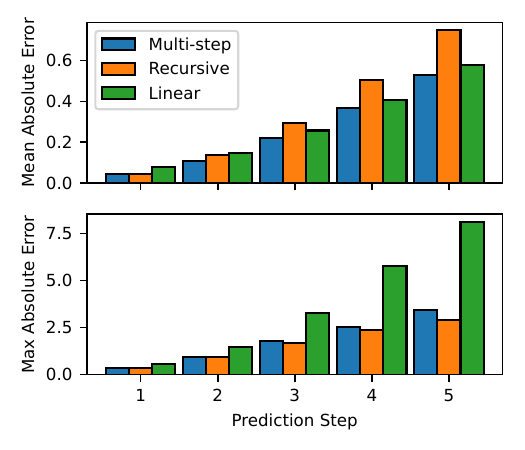}
    \caption{Mean and maximum absolute errors for $i = 1, \ldots, 5$-step ahead predictions for the multi-step DCNN, for recursive predictions using a single-step DCNN, and for the predictions of the linear ARI model}
    \label{fig:beta_prediction_errors}
\end{figure}

We first compare: (i) the multi-step predictions of the refined model, (ii) the recursive application of the single-step ahead predictor ($f_{1}(z_{k}, u_{k})$ in (\ref{eq:nn_dynamics})), and (iii) a linear model of beta oscillations (which is used in the linear MPC algorithm described in Section \ref{ssec:control_comparison}). For these tests, the predictors were trained and tested using patient activity with no  stimulation effect. As expected, prediction errors increase as the prediction horizon grows. Qualitatively, the multi-step predictor is less smooth than either of the recursive predictors (Figure~\ref{fig:predicted_beta_trajectories}) but is on average closer to the true beta activity. Figure~\ref{fig:beta_prediction_errors} shows that the multi-step predictor acheives the smallest mean absolute error across all prediction steps, but the maximum error is lower for the recursive predictor for $3$ or more prediction steps.

\subsection{Comparison with Alternative Control Strategies} \label{ssec:control_comparison}
This section compares the performance of DCNN TMPC with two control algorithms that have been tested \textit{in vivo}, and with the linear MPC strategy proposed in~\cite{steffen:2024}. We consider 50 simulations, each consisting of 15 seconds of patient data selected randomly from the test set, with the parameters of~(\ref{eqn:stim_dynamics}) varying stochastically as described in Section~\ref{ssec:stim_train}. 

The on-off controller, $u_{\mathrm{on-off}, k}$, increases the level of stimulation by the maximum increment $\Delta u_{max}$, up to the maximum value $u_{max}$ whenever the measured beta activity exceeds the threshold $\beta_0$, and reduces stimulation otherwise
\begin{equation}\label{eqn:onoff}
    u_{\mathrm{on-off}, k} = 
    \begin{cases}
        \mathrm{min}(u_{\mathrm{on-off}, k-1} + \Delta u_{max},  u_{max}) & \text{if}\ y_{k} > \beta_{0} \\
        \mathrm{max}(u_{\mathrm{on-off}, k-1} - \Delta u_{max},  0) & \text{otherwise.}
    \end{cases}
\end{equation}
The PI controller is implemented in difference form, $u_{\mathrm{PI}, k} = u_{\mathrm{PI}, k-1} + \Delta u_{\mathrm{PI}, k}$, using the error signal $e_{k} = [y(k T_s) - \beta_{0}]_{\geq0}$,
\begin{equation}
\label{eqn:pid}
\Delta u_{\mathrm{PI},k} = \  K_{P}\Delta e_{k} + K_{I}T_{s} e_{k} 
\end{equation}
where $K_{P}$ and $K_{I}$ are the controller gains, 
with $\Delta u$ and $u$ limited (via saturation) to $[-\Delta u_{max}, \Delta u_{max}]$ and $[0, u_{max}]$. 

The linear MPC controller is described in detail in \cite{steffen:2024}; here we describe only its main features. The approach uses a linearizing transformation of (\ref{eq:beta_response}), $\xi_{k} = \mathrm{ln}(y_{k})$, such that the effect of stimulation appears additively in the model, $\xi_{k} = \xi_{\beta, k} - \xi_{\eta, k}$, and a linear ARI model of nominal beta activity
\begin{equation}\label{eq:linear_beta}
    \Delta \xi_{\beta, k} = \sum_{i=1}^{n_{\beta}}\theta_{i}\Delta \xi_{\beta, k-1}
\end{equation} where $\Delta \xi_{\beta,k} = \xi_{\beta,k} - \xi_{\beta, k-1}$. The parameters $\theta_{1}, \ldots, \theta_{n_\beta}$ are identified from the patient data using a least squares approach over an identification period (during which there is no stimulation) immediately prior to initiating closed-loop control. The controller uses an augmented linear model
\begin{equation} \label{eqn:ss}
\begin{aligned}
    x_{k+1} &= A(\theta)x_{k} + Bu_{k} \\
    A &= \begin{bmatrix}
        A_{\eta} & {\mathbf 0}_{2 \times (n_{\beta}+1)} \\ 
        \mathbf{0}_{(n_{\beta}+1) \times 2} & A_{\beta}(\theta)
    \end{bmatrix} \\
    B &= \begin{bmatrix}
        B_{\eta} \\
        \mathbf{0}_{(n_{\beta}+1) \times 1}
    \end{bmatrix} \\
    \xi_{k}  &= \begin{bmatrix}
        0 & 1 & -1 & \mathbf{0}_{1 \times n_{\beta}}
    \end{bmatrix}  x_{k}
\end{aligned}
\end{equation}
where the matrices $A_{\eta}$ and $B_{\eta}$ are computed by discretizing the model (\ref{eqn:stim_dynamics}) with the average parameter values and a zero-order hold, and $A_{\beta}(\theta)$ is the linear nominal activity model (\ref{eq:linear_beta}) in state-space form. The controller solves the optimal control input for a quadratic cost index
\begin{equation}
    \underset{u}{\mathrm{argmin}} \sum_{i=0}^{N-1} \Bigl( [\xi_{k+i}-\xi_{0}]_{\geq0}^{2} + Ru_{k+i}^{2} \Bigr)
\end{equation} subject to (\ref{eqn:ss}), $\Delta u_{k+i} \in \Delta \mathbb{U}$ and $u_{k+i} \in \mathbb{U}$ for all $i \in \{ 0, \ldots, N-1\}$.

The DCNN TMPC algorithm was implemented in Python using the cvxpy library \cite{diamond2016cvxpy} to interface with the MOSEK solver. The maximum output constraint $y_{max}$ was chosen as the 95\textsuperscript{th} percentile of the patient's beta activity, while $y_{min}$ was chosen as the minimum recorded value. We used a prediction horizon of $N=5$ steps, and the disturbance set $\mathbb{W}_{i}$ for each $i\in\{1,\ldots,5\}$ was chosen as the 80\textsuperscript{th} percentile absolute prediction error during training. We note that this constraint cannot guarantee robust output constraint satisfaction. However, this is not problematic in practice because the safety of the DBS system is guaranteed by appropriately chosen constraints on the input $u$ and rate of change $\Delta u$, so violations of output constraints can be handled by progressive constraint softening until the problem once again become feasible. 

\begin{figure}
    \centering
    \includegraphics[width=0.85\linewidth,trim={0 4mm 0 3mm},clip]{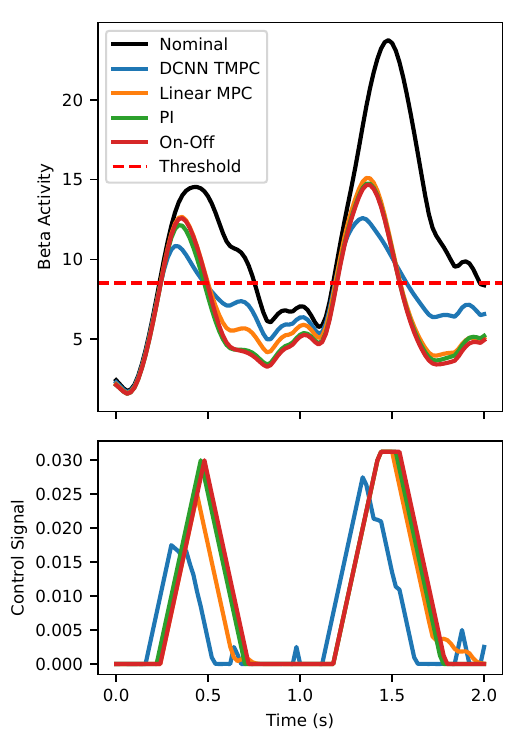}
    \caption{Comparison of trajectories and control input sequences for DCNN TMPC (using the refined model), linear MPC, PI and on-off controllers. Here `Nominal' indicates the raw beta activity with no applied stimulation}
    \label{fig:trajectories}
\end{figure}

\begin{figure}
    \centering
    \includegraphics[width=0.85\linewidth,trim={0 4mm 0 3mm},clip]{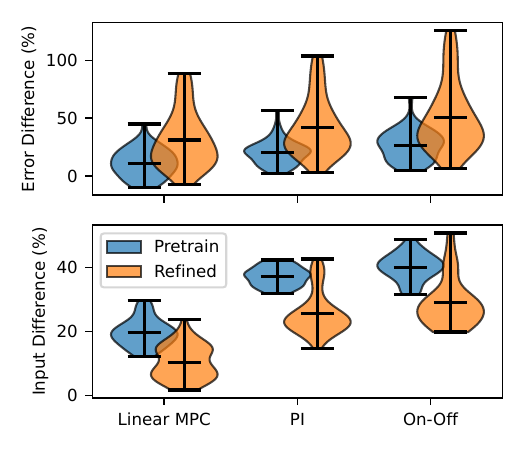}
    \caption{Violin plot of percentage improvement in beta suppression error and control input for the DCNN TMPC algorithm with pretrained and refined model, relative to linear MPC, PI and on-off controllers for the 50 simulations. Bars show the mean, min, and max of the data}
    \label{fig:controller_metrics}
\end{figure}

\begin{figure}
    \centering
    \includegraphics[width=0.85\linewidth,trim={0 4mm 0 3mm},clip]{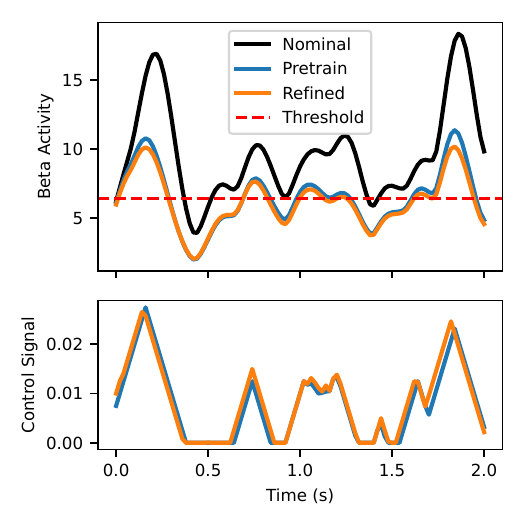}
    \caption{Comparison of trajectories and control sequence for DCNN TMPC using a model with pre-training only, versus the refined model}
    \label{fig:refined_vs_pretrained_trajectories}
\end{figure}

Figure \ref{fig:controller_metrics} shows that DCNN TMPC provides significant improvements in performance compared to the linear MPC, PI and on-off controllers. With the refined model, DCNN TMPC beta suppression error is on average 30-50\% higher for the alternative controllers. In addition, DCNN TMPC with only the pre-trained model outperforms linear MPC by at least 10\% on average, and outperforms PI and on-off controllers by at least 20\%. There is a smaller variation in the applied control input, however in all simulations, DCNN TMPC is more than 5\% more energy-efficient than linear MPC, and 20\% better than PI and on-off control. This highlights the benefits of using a more expressive nonlinear model for predicting the future system behavior. Moreover, the control input sequences in Fig.~\ref{fig:trajectories} show that DCNN TMPC 
responds more quickly than the alternative controllers to bursts in beta activity 
and reduces stimulation faster in response to decreases in beta activity. The beta activity plots also show that the proposed controller uses less unnecessary stimulation; this can be seen from the fact that the DCNN TMPC modulated trajectory is much closer to the nominal activity when it is below the threshold. 

Comparing performance with pre-trained and refined models reveals significant robustness to differences in beta activity across different patients. Figure~\ref{fig:controller_metrics} shows that even with a model that was not trained on data from the test patient, DCNN TMPC is able to perform significantly better than the alternative controllers. 
Remarkably this observation even applies when DCNN TMPC with a pre-trained model is compared with linear MPC, which uses a model trained on data from the test patient.
The comparison of the trajectories of DCNN TMPC with pre-training and refined models 
in Fig.~\ref{fig:refined_vs_pretrained_trajectories} shows very similar trajectories, with the refined model applying slightly more control action, and generally acting faster, which results in better suppression of pathological beta activity. Generally it appears that the model with pre-training alone results in a less aggressive controller.

\section{Concluding Remarks and Further Work} \label{sec:conclusions}
This paper demonstrates that
DCNN TMPC with a multi-step predictor outperforms existing control algorithms for CLDBS in simulations using actual patient data. We demonstrate that the multi-step NN predictor performs better than recursive predictions, and that the neural network-based controller generalizes to differences in beta activity across different patients.

We believe our approach provides a readily generalizable framework for CLDBS, and indeed many other closed-loop neuromodulation techniques, as it makes few assumptions about the dynamic model underlying the response of the relevant bio-marker. This approach will remain useful even if a different biomarker (or set of biomarkers) is discovered to be a better indicator of disease state.

A limitation of the multi-step predictor approach is that it is unable to ensure recursive feasibility since there is no guarantee that the errors in multi-step predictions will be consistent for predictions made at different times. We aim to investigate methods of bounding prediction errors such that we can develop either robust or stochastic feasibility certificates. In addition, we plan to validate controller performance in biophysical models and \textit{in vivo} with PD patients. 

\bibliographystyle{ieeetr}
\bibliography{dcnndbs}

\end{document}